\documentclass[a4paper,10pt]{article}
\usepackage[margin=1.0in]{geometry}
\usepackage{amsmath,amsthm,amssymb,mathtools,tikz-cd,dsfont}
\usepackage[colorlinks,citecolor=red,urlcolor=blue,bookmarks=true,hypertexnames=true]{hyperref}
\hypersetup{pdftitle={Derived Equivalences for Terminal Threefold Flops via Local Index-One Covering Models},pdfauthor={Hao Xingbang}}
\newcommand{\X}{\mathrm{X}}
\newcommand{\W}{\mathrm{W}}
\newcommand{\V}{\mathrm{V}}
\newcommand{\Pj}{\mathrm{P}}
\newcommand{\OO}{\mathcal O}
\newcommand{\Db}{\mathrm{D}^b}
\newcommand{\id}{\mathrm{id}}
\newcommand{\Cone}{\operatorname{Cone}}
\newcommand{\Spec}{\operatorname{Spec}}
\newtheorem{lemma}{Lemma}
\newtheorem{theorem}{Theorem}
\theoremstyle{remark}
\newtheorem{remark}{Remark}
\theoremstyle{plain}
\allowdisplaybreaks[2]
\begin{document}
\begin{center}
\fontsize{13pt}{10pt}\selectfont
\textsc{\textbf{DERIVED EQUIVALENCES FOR TERMINAL THREEFOLD FLOPS\\[3pt]
VIA LOCAL INDEX-ONE COVERING MODELS}}
\end{center}
\vspace{0.1cm}
\begin{center}
\fontsize{10pt}{8pt}\selectfont
\textsc{HAO XINGBANG}
\end{center}
\vspace{0.2cm}
\begin{center}
\fontsize{12pt}{10pt}\selectfont
\textsc{Abstract}
\end{center}
\noindent
Inspired by the quadratic involutions of Corti--Pukhlikov--Reid
\cite{CPR} and Gavin Brown's classification of hypersurface flips via
toric GIT \cite{Brown}, we apply similar symmetries to the local covering
models of terminal threefold flops. We recover the derived equivalence
of their canonical stacks by pullback and pushforward along the fibre
product over the canonical stack of the contraction target, revisiting
the local picture in Kawamata's work \cite{Kawamata}.

\medskip
\noindent\textbf{AI disclosure.}
The ideas in this note originate with the author. AI tools were used
to help check and refine the proofs and to assist with the writing.

\smallskip
\noindent\textbf{Acknowledgements.}
The author thanks Hsin-Ku Chen for discussions and apologizes to him regarding their
joint project using his elegant factorization \cite{HKChen}, an approach that remains
highly valuable.

\medskip
All schemes and stacks are over $\mathds{C}$; $\Db$ denotes the bounded derived
category of coherent sheaves. Pullback and pushforward functors are
derived. Derived tensor products
are written $\otimes^{\mathbf L}$.
For local models, we follow Koll\'ar's classification \cite{KollarSurvey}.

\begin{lemma}[Quadratic model]\label{lem:quadratic}
Let $f:\X\to\W\leftarrow\X^+:g$ be a flop germ of terminal
Gorenstein threefolds, allowing reducible exceptional fibres, where
\[
\W=\Spec\frac{\mathds{C}[[z_1,z_2,z_3,s]]}{(s^2-t)},
\qquad
\Pj=\Spec\mathds{C}[[z_1,z_2,z_3]],
\qquad t\in\mathds{C}[[z_1,z_2,z_3]].
\]
Set $\V=\X\times_\W\X^+$, and let
$\operatorname{pr}:\W\to\Pj$ forget $s$; the two maps from the larger
product to $\W$ agree precisely on $\V$.
\[
\begin{tikzcd}[column sep=3.2em,row sep=1.1em]
&\V\arrow[ddl,bend right=30,"\phi"']
\arrow[d,hook]\arrow[ddr,bend left=30,"\varphi"]&\\
&\X\times_\Pj\X^+
\arrow[dl,"\operatorname{pr}_1"']
\arrow[dr,"\operatorname{pr}_2"]&\\
\X\arrow[dr,"f"]
\arrow[ddr,bend right=35,"\eta"']
&&\X^+\arrow[dl,"g"']
\arrow[ddl,bend left=35,"\eta^+"]\\
&\W\arrow[d,"\operatorname{pr}"]&\\
&\Pj&
\end{tikzcd}
\qquad
\eta=\operatorname{pr}\circ f,\quad
\eta^+=\operatorname{pr}\circ g.
\]
Then
\begin{equation}\label{eq:equivalence}
\varphi_*\phi^*:\Db(\X)\xrightarrow{\sim}\Db(\X^+),
\qquad
(\varphi_*\phi^*)^{-1}=\phi_*\varphi^!.
\end{equation}
\end{lemma}

\begin{proof}
\begin{enumerate}
\item \emph{Tor base change.}
Terminal Gorenstein threefolds are locally hypersurfaces
\cite{Reid}. On an open $\X^\circ\subset\X$, choose a regular immersion
\[
i:\X^\circ\hookrightarrow U\subset\mathbb A^r_\Pj\longrightarrow\Pj,
\qquad
\mathcal I_{\X^\circ/U}=(h_1,\ldots,h_r),
\]
with $U$ open. The Koszul resolution on $U$ is
\[
0\longrightarrow\bigwedge^r\OO_U^{\oplus r}
\longrightarrow\cdots\longrightarrow\OO_U^{\oplus r}
\xrightarrow{(h_1,\ldots,h_r)}\OO_U
\longrightarrow i_*\OO_{\X^\circ}\longrightarrow0.
\]
Its terms $\bigwedge^j\OO_U^{\oplus r}$ are $\Pj$-flat since
$U\to\Pj$ is smooth. Base change by $\X^+\to\Pj$ replaces it by
\[
0\longrightarrow\bigwedge^r\OO_{U\times_\Pj\X^+}^{\oplus r}
\longrightarrow\cdots\longrightarrow\OO_{U\times_\Pj\X^+}^{\oplus r}
\xrightarrow{(h_1,\ldots,h_r)}\OO_{U\times_\Pj\X^+}.
\]
We should check that this complex is exact in negative degrees.
At a closed point $\xi$ of $\X\times_\Pj\X^+$ lying in this chart,
$B=\OO_{U\times_\Pj\X^+,\xi}$ is Gorenstein of dimension $3+r$,
since $U\times_\Pj\X^+\to\X^+$ is smooth of relative dimension $r$.
The upper bound below holds because $B/(h_1,\ldots,h_r)$ is a local
ring of $\X^\circ\times_\Pj\X^+$: components not lying over $0$ are
generically finite over their images in $\Pj$, while the fibre over $0$
is a product of two fibres of dimension at most one, hence at most two.
Together with Krull's height theorem, this gives
\[
3=\dim B-r\le\dim B/(h_1,\ldots,h_r)\le3.
\]
Since $B$ is Cohen--Macaulay, this dimension equality makes
$(h_1,\ldots,h_r)$ a regular sequence \cite[\href{https://stacks.math.columbia.edu/tag/02JN}{Tag~02JN}]{Stacks},
so the base-changed Koszul complex resolves $B/(h_1,\ldots,h_r)$.
This proves the asserted exactness at closed points. Since
$\X\times_\Pj\X^+$ is proper over $\Pj$, it follows that
\begin{equation}\label{eq:tor}
\mathrm H^{-i}\!\left(\OO_\X\otimes^{\mathbf L}_{\OO_\Pj}\OO_{\X^+}\right)=0
\quad(i>0),
\qquad
\eta^{+*}\eta_*\simeq
\operatorname{pr}_{2*}\operatorname{pr}_1^*.
\end{equation}
The second identity is Tor-independent base change
\cite[\href{https://stacks.math.columbia.edu/tag/08IB}{Tag~08IB}]{Stacks}. The same argument applies to $\X\times_\Pj\X$;
both products are Gorenstein.

\item \emph{Root exchange and the kernel sequence.}
The involution $\tau(s)=-s$ acts by $-1$ on $\mathrm{Cl}(\W)$:
for any Weil divisor $\mathrm{D}$,
\[
\mathrm{D}+\tau^*\mathrm{D}\sim\operatorname{pr}^*\operatorname{pr}_*\mathrm{D},
\]
and the right side is principal because $\Pj$ is factorial.
Smallness identifies the divisorial section algebras on $\W$ with
those on the two models. Relative Proj and the opposite ampleness
assumption therefore give \cite[Proposition 2.2, Example 2.3]{Kollar}
\[
\sigma_\X:\X^+\xrightarrow{\sim}\X,
\qquad f\sigma_\X=\tau g,
\qquad \eta^+=\eta\sigma_\X.
\]

The fibre product $\X\times_\Pj\X^+$ is Cohen--Macaulay and generically
reduced, with precisely the matching component and the graph
$\Gamma_{\sigma_\X}$ as its generic components. The dimension argument
above excludes additional vertical components, so it is reduced.
The ideal of $\V$ is $(s-s^+)$, and
$(s-s^+)(s+s^+)=0$. Reducedness shows that $s-s^+$ annihilates the
ideal of $\Gamma_{\sigma_\X}$. Its restriction to this integral graph
is $-2s^+$, a nonzero divisor. Thus multiplication gives the exact sequence
\begin{equation}\label{eq:kernel}
0\longrightarrow\OO_{\Gamma_{\sigma_\X}}
\xrightarrow{\ s-s^+\ }
\OO_{\X\times_\Pj\X^+}
\longrightarrow\OO_\V\longrightarrow0.
\end{equation}
Dualizing the diagonal restriction gives the unit of $\eta_*\dashv\eta^*$.
On the finite locus, the functional $a+bs\mapsto b$ gives
$1\mapsto1\otimes s+s\otimes1$; root exchange gives the multiplier
$s-s^+$. Reducedness determines maps from the graph kernel by their
generic restrictions, so this is the first map in \eqref{eq:kernel} globally.
Taking Fourier--Mukai transforms in \eqref{eq:kernel} gives
\begin{equation}\label{eq:unit-cone}
\varphi_*\phi^*
\simeq\sigma_\X^*\Cone\bigl(
\id\longrightarrow \eta^*\eta_*\bigr),
\end{equation}
where the arrow is the unit of $\eta_*\dashv \eta^*$.

\item \emph{The two adjunctions.}
Since $\Pj$ is regular and $\eta$ is proper, $\eta^*$ and $\eta_*$
preserve $\Db$, and hence so do both cone functors.
Crepant duality gives $\eta^!\OO_\Pj\simeq\OO_\X$
\cite[\href{https://stacks.math.columbia.edu/tag/0B6U}{Tag~0B6U}]{Stacks},
and rational singularities give $f_*\OO_\X\simeq\OO_\W$.
The projection-formula map is an isomorphism on $\OO_\Pj$, which
generates $\Db(\Pj)=\operatorname{Perf}(\Pj)$, so the resulting
identities are natural on the derived categories:
\[
\eta^*\dashv \eta_*\dashv \eta^*,
\qquad \eta_*\eta^*\simeq\id\oplus\id.
\]
Write $u$ for the unit of $\eta_*\dashv\eta^*$ and $\varepsilon$ for
the counit of $\eta^*\dashv\eta_*$. The basis $1,s$ splits
$(\eta^*\eta_*)^2\simeq(\eta^*\eta_*)^{\oplus2}$. At the level of
Fourier--Mukai kernels, one composite is therefore
\[
\begin{gathered}
\Cone\bigl(\eta^*\eta_*\xrightarrow{\varepsilon}\id\bigr)
\circ\Cone\bigl(\id\xrightarrow{u}\eta^*\eta_*\bigr)[-1]
\\
\simeq\left[
\eta^*\eta_*
\xrightarrow{\left(\begin{smallmatrix}\eta^*\eta_*(s)\\1\\\varepsilon\end{smallmatrix}\right)}
(\eta^*\eta_*)^{\oplus2}\oplus\id
\xrightarrow{\left(\begin{smallmatrix}1&s&-u\end{smallmatrix}\right)}
\eta^*\eta_*
\right]\simeq\id.
\end{gathered}
\]
The three terms lie in degrees $-1,0,1$; $s$ acts on the target factor
and $\eta^*\eta_*(s)$ on the source factor.
Since $u\varepsilon=s+\eta^*\eta_*(s)$, elimination cancels the two
identity blocks and leaves $\id$ in degree zero. The reverse composite
has the same cancellation with the two $s$-actions interchanged.
Taking right adjoints in \eqref{eq:unit-cone} before restricting to $\Db$ gives
$\phi_*\varphi^!\simeq
\Cone(\eta^*\eta_*\xrightarrow{\varepsilon}\id)[-1]\circ\sigma_{\X*}$,
since the adjoint of $u$ is $\varepsilon$.
Thus the local pull--push functor and its inverse preserve $\Db$,
and \eqref{eq:unit-cone} proves \eqref{eq:equivalence}.

\end{enumerate}
\end{proof}

\begin{lemma}[Brown--Mallinson--Chang model]\label{lem:diagonal}
Let $\X\dashrightarrow\X^+$ be the terminal flop obtained by exchanging
$g_1$ and $g_2$ over
\[
\W=\{xy-g_1(z,u)g_2(z,u)=0\},
\]
where $(g_1,g_2)$ is a regular sequence. Set $\V=\X\times_\W\X^+$,
with projections $\phi:\V\to\X$ and $\varphi:\V\to\X^+$. Then
\[
\varphi_*\phi^*:\Db(\X)\xrightarrow{\sim}\Db(\X^+),
\qquad
(\varphi_*\phi^*)^{-1}=\phi_*\varphi^!.
\]
\end{lemma}
\begin{proof}
Koll\'ar's construction \cite[(2.2.2.2)]{KollarSurvey} blows up
$(x,g_2)$ and $(x,g_1)$, respectively. Their graded GIT presentation is
\cite{Brown,Chang,Mallinson}
\begin{equation}\label{eq:diagonal-cox}
\frac{\mathds{C}[[z,u]][x_1,x_2,y_1,y_2]}
{(x_1y_1-g_1,\ x_2y_2-g_2)},
\qquad
\deg(x_1,x_2,y_1,y_2)=(1,1,-1,-1).
\end{equation}
The two chambers are $\mathbf x\ne0$ and $\mathbf y\ne0$;
$x=x_1y_2$ and $y=x_2y_1$.
As in \cite{HaoIIa}, blow up the ambient centre $(y_1,y_2)$ with
weights $(1,1)$. The two ambient projections are
\[
\bar\phi:(\mathbf x,\mathbf y,t)\longmapsto(\mathbf x,t\mathbf y),
\qquad
\bar\varphi:(\mathbf x,\mathbf y,t)\longmapsto(t\mathbf x,\mathbf y).
\]
Both total transforms are therefore
\begin{equation}\label{eq:diagonal-roof}
\frac{\mathds{C}[[z,u]][x_1,x_2,y_1,y_2,t]}
{(tx_1y_1-g_1,\ tx_2y_2-g_2)},
\qquad
\deg x_i=(1,0),\quad\deg y_i=(0,1),\quad\deg t=(-1,-1),
\end{equation}
modulo $\mathbb G_m^2$, on $\mathbf x\ne0$, $\mathbf y\ne0$.
Modulo $t$, the equations are the regular sequence $-g_1,-g_2$; hence
\[
(tx_1y_1-g_1,\ tx_2y_2-g_2):t^\infty
=(tx_1y_1-g_1,\ tx_2y_2-g_2).
\]
Thus strict and total transforms agree, and the ambient base changes
are Tor-independent. The two blow-ups specify the column and row lines of
\[
\begin{pmatrix}g_1&x\\y&g_2\end{pmatrix}
=t\binom{x_1}{x_2}(y_1\ y_2),
\]
which identifies the common transform with $\X\times_\W\X^+$,
including its scheme structure. The Type~IV standard-flop calculation
in \cite{HaoIIa}, with $1+1-1-1=0$, gives the equivalence.
\end{proof}

\medskip
\begin{remark}[A trivial formal equivalence]
For the Gorenstein local covering models, with exceptional fibres $C$
and $C^+$, Koll\'ar's involution \cite{Kollar} gives
\[
\widehat{\X}_C\simeq\widehat{\X}^+_{C^+},
\qquad
\Db(\widehat{\X}_C)\simeq\Db(\widehat{\X}^+_{C^+}).
\]
The involution need not extend to the global models.
\end{remark}

\begin{theorem}[Kawamata {\cite[Theorem~6.5]{Kawamata}}]\label{thm:terminal}
Let $f:\X\to\W\leftarrow\X^+:g$ be a flop of terminal threefolds,
with $f$ and $g$ projective. For their canonical stacks \cite{HaoIIa,HaoIIb}, let
\[
\widetilde\X\xleftarrow{\ \phi\ }
\widetilde\V=\widetilde\X\times_{\widetilde\W}\widetilde\X^+
\xrightarrow{\ \varphi\ }\widetilde\X^+.
\]
Then
\[
\varphi_*\phi^*:\Db(\widetilde\X)
\xrightarrow{\sim}\Db(\widetilde\X^+),
\qquad
(\varphi_*\phi^*)^{-1}=\phi_*\varphi^!.
\]
\end{theorem}
\begin{proof}
\begin{enumerate}
\item \emph{Local models.}
Pass to the index-one cover of the completed contraction base.
Small crepantness gives terminal Gorenstein covering models over an
isolated cDV singularity \cite[Corollary 2.2.21, Theorem 2.3.4]{Kawakita}.
Its equation has the quadratic form of Lemma~\ref{lem:quadratic}
\cite[Proposition 2.3.3]{Kawakita}, which therefore applies to both
irreducible and reducible fibres.
For irreducible fibres, Koll\'ar's classification
\cite[Theorem 2.2.2]{KollarSurvey} gives the explicit models of
Lemmas~\ref{lem:quadratic} and~\ref{lem:diagonal}.
The local reduction \cite[Lemmas 5.1.4--5.1.5]{Kawakita}
also provides a factorization viewpoint.

\item \emph{Descent.}
Locally, crepant smallness identifies the correspondence with
\[
\widetilde\V=
\bigl[(\X^\#\times_{\W^\#}\X^{+\#})/\mu_r\bigr]
\]
for the diagonal action \cite[proof of Theorem 6.5]{Kawamata}.
The index-one covering maps are \'etale in codimension one, so they
introduce no ramification divisor.
The kernel and its adjunction maps are equivariant, so the local
equivalences descend. The cone descriptions give uniform cohomological
bounds for the pull--push functor and its inverse; these bounds descend
through faithfully flat completion and the quotient atlases. Away from
the exceptional locus the correspondence is a graph. Thus the global
adjunction $\varphi_*\phi^*\dashv\phi_*\varphi^!$ restricts to $\Db$.
Its unit and counit are isomorphisms after completion at the exceptional
values and away from them. Faithful flatness proves they are
isomorphisms globally. The admissibility and spanning-class argument
of \cite{HaoIIa} thus gives the global equivalence for this same
pull--push functor.
\end{enumerate}
\end{proof}

\medskip
\noindent\textbf{Remark.}
The smooth and Gorenstein cases were studied by Bridgeland
\cite{Bridgeland}, Chen \cite{Chen} and Van den Bergh \cite{VdB},
with the terminal setting treated by Kawamata \cite{Kawamata} through
canonical stacks. Perverse hearts on unbounded quasi-coherent categories,
``Serre'' spectral sequences, and tilting theory on singular varieties
do not make boundedness automatic; compare Bodzenta--Bondal
\cite[Remark~2.3]{BBFlops} and \cite{BBTilting}.

\bigskip
{\footnotesize
\noindent HAO XINGBANG, \textsc{Waseda University, Ookubo, Shinjuku-ku,
Tokyo, 169-8555, Japan.}\par\nopagebreak
\noindent\textit{E-mail address}: \texttt{hao@fuji.waseda.jp}\par
}
\end{document}